\documentclass[11pt,twoside]{article}
\usepackage{graphicx} % Required for inserting images
\usepackage{ulem}
\usepackage{amssymb}
\usepackage[mathscr]{eucal}
\usepackage{eufrak}
\usepackage{amsmath,amsthm}
\usepackage{physics}
\usepackage{mathrsfs}
\usepackage{lgreek}
\usepackage[shortlabels]{enumitem}
\usepackage[colorlinks=true,
   urlcolor=blue,           filecolor=green,      
   citecolor=green,      
   linkcolor=red,           bookmarks=true,
  unicode,
   plainpages=false,   ]{hyperref}
\usepackage{color}
\usepackage{comment}
\font\bg=cmbx10 scaled 1300

\definecolor{royalblue}{rgb}{0,0,0.128}
\def\bp{\begin{proof}}
\def\ep{\end{proof}}
\def\n{\nabla}

\def\dm{|\hskip-0.05cm|}

\def\OO{\Omega}
\def\displ{\displaystyle}

\def\VS{\vspace{6pt}\\\displ }

\def\R{\Bbb R}
\def\N{\Bbb N}
\def\à{à}

\def\vep{\varepsilon}

\def\be{\begin{equation}}
\def\ba{\begin{array}}
\def\ea{\end{array}}
\def\ee{\end{equation}}
\def\vs1{\vspace{1ex}}

\def\Ã©{\'{e}}
\def\Ãš{\`{e}}

\newtheorem{lemma}
{\bf Lemma} 
\font\sc=cmcsc10
\newtheorem{ass}
{\bf Assumption} 
\newtheorem{defi}
{\bf Definition} 
\newtheorem{tho}
{\bf Theorem} 
\newtheorem{rem}
{\sc Remark}

\title{\bg Navier-Stokes flow past a translating body: \\ attainability of the steady state as limit of an $L^p$ weak solution}
\author{\sc  Filippo Palma
\thanks{Dipartimento di Matematica e Fisica,  
Universit\`{a} degli 
Studi della Campania
``L. Vanvitelli'', via Vivaldi 43, 81100 \null\hskip0.55cmCaserta,
 Italy.\newline\null\hskip0.55cm
 filippo.palma@unicampania.it
\newline\null\hskip0.55cm The  research activity  is performed under the
auspices of   GNFM-INdAM. }}
\begin{document}
\markboth{\footnotesize\rm    F. Palma} {\footnotesize\rm Navier-Stokes flow past a translating body\dots
}
\maketitle
\begin{abstract}
    Consider a viscous, incompressible, Newtonian fluid that fills the region of space exterior to a rigid body. Let the fluid be activated by a perturbation $u_0\in L^p(\OO)$, $p\in (2,3)$ and let the body start to translate with a prescribed velocity $-h(t)u_\infty$, $t>0$, $u_\infty\in\R^3\setminus \{0\}$. Suppose that the body reaches the stationary regime after some finite time. We show that, under suitable assumptions on $u_\infty$, a steady state is always reachable by the fluid, asymptotically in time, as limit of a suitable weak solution to the Navier-Stokes flow in suitable metrics. Our result seems to be a further contribution related to the {\it Starting Problem}, posed by R. Finn in 1965.
\end{abstract}
\section{Introduction}
Let us consider a rigid body $\mathcal  B$ that translates in a viscous, incompressible, Newtonian fluid. By the symbol $\mathcal B$ we also identify the region of space occupied by the body, assuming it to be compact and with a smooth boundary $\partial\mathcal  B$. We assume that the fluid fills the entire region of space exterior to the body, and we follow the evolution of the fluid flow in a frame attached to the body. \par Setting $\OO=\R^3\setminus \mathcal B$, the governing equations of the fluid motion, in the aforementioned frame, read as follows:
\begin{equation}\label{NS-Ose}
    \begin{cases}
        u_t-\Delta  u+u\cdot\n u+hu_\infty\cdot\n u+\n\pi_u=0\,,\,\,\text{in }(0,T)\times\OO\,,\\ \hskip5.3cm\n\cdot u=0\,,\,\,\text{in }(0,T)\times\OO\,,\\ \hskip3.75cmu(t,x)=-h(t)u_\infty\,,\,\,\text{on }(0,T)\times\partial\OO\,, \\ \hskip2.8cm\displ\lim_{|x|\to\infty }u(t,x)=0\,,\\ \hskip2.8cm(u,\varphi)=(u_0,\varphi)\,,\,\,\text{for }t=0\,,\,\,\text{for all }\varphi\in\mathscr C_0(\OO)\,.
    \end{cases}
\end{equation}
In problem \eqref{NS-Ose}, $u:(0,T)\times\OO\to \R^3$ is the kinetic field of the fluid, $\pi_u:(0,T)\times\OO\to\R$ is the pressure field. Concerning the function $h$, we make  the following assumption:
\begin{ass}\label{hfunct}
{\sl
\[
\ba{c}
h\in C^{1,\alpha}([0,\infty))\,,\,\,\alpha\in (0,1)\,,\,\, h(0)=0\,,\VS \text{there exists }T_0>0 \,\text{such that }h(t)=1\,\text{for all }t\geq T_0\,.
\ea
\]}
\end{ass}
The aim of this work is to show that, under Assumption\,\ref{hfunct} and provided $|u_\infty|$ is sufficiently small, for all data $u_0\in L^p(\OO)$, $p\in(2,3)$, problem \eqref{NS-Ose} admits a weak solution $u$ (see Definition\,\ref{def: WS-Ose}) which converges, asymptotically in time and in suitable metrics, to a stationary regime. \par The investigation of the asymptotic properties is justified by the fact that the weak solution we obtain is regular a. e. in time and, in particular, it is regular for times $t\geq \theta$, $\theta\geq 0$ suitable.\par In the framework of {\it physically reasonable} stationary solutions introduced by Finn, see \cite{Finn-sta} and further references therein, the stationary regime reached by our solution asymptotically in time is given by the unique solution to
\begin{equation}\label{NS-Ose-staz}
    \begin{cases}
        -\Delta U+U\cdot\n U+u_\infty\cdot\n U+\n\pi_U=0\,,\,\,\text{in }\OO\,,\\ \hskip4.8cm\n\cdot U=0\,,\,\,\text{in }\OO\,,\\\hskip4.9cm U|_{\partial\OO}=-u_\infty\,,\\\hskip3.96cm \displ\lim_{|x|\to\infty} U(x)=0\,,
    \end{cases}
\end{equation}
provided that a smallness assumption on $u_\infty$ is satisfied. \par
We give the following notion of a weak solution to problem \eqref{NS-Ose} (for details related to the notation, see section\,\ref{sec: prel}).
\begin{defi}\label{def: WS-Ose}
    {\sl  Let $u_0\in L^p(\OO)$, $p\in (2,3]$. A field $u:(0,\infty)\times\OO\to \R^3$ is said a weak solution to problem \eqref{NS-Ose} if, for all $T>0$,
    \begin{description}
        \item[a)] $u\in L^{2,\infty}(0,T;W^{1,2}_{loc}(\OO))\cap L^{4,\infty}(0,T;L^3(\OO))$ and $(u(t),\n\psi)=0$, for all $\psi\in C_0^\infty(\OO)$ and $t\in (0,T)$;
        \item[b)] $$\ba{c}\displ \int_0^T \big[(u,\phi_t)-(\n u,\n \phi)+((u+u_\infty)\cdot\n \phi,u)\big]\,d\tau+(u_0,\phi(0))=0\,,\VS \text{for all }\phi\in C_0^\infty([0,T)\times\OO)\,,\,\,\text{with }\n\cdot\phi=0\,\,\text{for all }t\in [0,T)\,;\ea$$
        \item[c)] $\gamma (u(t))=-h(t)u_\infty$ in $H^\frac12(\partial\OO)$ for a. a. $t\in (0,T)$ and $\displ\lim_{t\to 0}(u(t),\varphi)=(u_0,\varphi)$, for all $\varphi\in \mathscr C_0(\OO)$.
    \end{description}}
\end{defi}
Our main result is the following.
\begin{tho}\label{thm: MR}
    {\sl Let $p\in (2,3)$. There exists $\delta>0$ such that: if $u_\infty\in\R^3\setminus \{0\}$ satisfies $|u_\infty|\leq \delta$, then, for all $u_0\in L^p(\OO)$, and for any function $h(t)$ satisfying Assumption\,\ref{hfunct}, there exists a weak solution $u$ to problem \eqref{NS-Ose} verifying the asymptotic property
    \be\label{As-Finn}
\dm u(t)-U\dm_q = O(t^{-\frac32(\frac13-\frac1q)})\,,\,\,\text{as }t\to\infty\,,\,\,\text{for all }3\leq q\leq \infty\,,
\ee
where $U$ is the unique physically reasonable solution to problem \eqref{NS-Ose-staz}.}
\end{tho}
\begin{rem}
    {\rm \begin{description}
        \item[(i)] In the particular case in which the body is at rest, hence $u_\infty=0$, the existence of a weak solution for initial data in $L^p(\OO)$, $p\in (2,3)$, was established in the recent paper \cite{Pal}. In accord with the quoted paper, we are able to obtain for the weak solution to problem \eqref{NS-Ose} a Leray-type \textit{Th\'eor\`eme de Structure}. This allows us to get the asymptotic property \eqref{As-Finn}.  
        \item[(ii)] The case $p=3$ was investigated in \cite{His-Mar}. We claim that our argument lines, suitably modified, can be adapted also to this particular case, leading to an alternative proof of the result therein. We will give more details in the following of the paper.
        \item[(iii)] In the statement of the theorem, we have implicitly assumed the well-posedness of the stationary problem \eqref{NS-Ose-staz}. In Lemma\,\ref{le: WP-staz} we will recall this result. 
        \item[(iv)] Notice that Theorem\,\ref{thm: MR} is, first of all, an existence result for a weak solution to problem \eqref{NS-Ose}. Hence, we are not claiming that the steady state is always reached by the fluid, but it is always reachable through a particular path. This is related to the lack of uniqueness of a weak solution to the Navier-Stokes equations.
        \item[(v)] We point out that the rate of asymptotic convergence to the stationary regime for our weak solution is different from that established in \cite{Pal} in the case $U=0$. However, this seems to be in accord with the results available for perturbations of a Navier-Stokes fluid from a steady state, see \cite{Mar-AB}.
    \end{description}}
\end{rem}
\subsection{Literature overview}
In the paper \cite{Finn-AMS}, R. Finn raised a very interesting question related to the attainability of a steady state for a fluid that interacts with a moving obstacle. Heuristically, the problem posed by Finn can be summarized as follows:
\begin{center}
    {\it Suppose that a body, which is originally at rest, is accelerated and starts to move in a fluid that occupies the region of space exterior to it. Let the body reach a stationary regime after some finite time. In a frame attached to the body, does the fluid reach a steady state asymptotically in time?}
\end{center}
The aforementioned problem is known in the literature as the Finn {\it Starting Problem} and it has been the subject of intensive research, as it remained essentially an open problem for more than thirty years. \par If the fluid is initially at rest and activated only by the motion of the body, a resolution of the problem was found by Galdi, Heywood and Shibata in \cite{GHS}, who proved the unique solvability of the system related to the fluid motion and the asymptotic convergence in time of the kinetic field of the fluid to a stationary regime, giving also a suitable convergence rate. The convergence is meant in the $L^\infty$-norm. \par The proof given in the quoted paper relies on the properties of the Oseen semigroup found in \cite{Kob-Shi} that allowed to find a global solution to a special exterior nonstationary problem for the Navier-Stokes equations, see \cite{Eno-Shi, Shi}. The peculiarity of this exterior problem consists in the fact that the fluid is required to tend to a prescribed nonzero constant velocity at spatial infinity. \par An improved estimate of the convergence rate to the stationary regime was given by Takahashi in \cite{Tak}. In the same paper, the author also discusses the Stating Problem in higher dimensions. \par  A natural question one might ask is the following: {\it what is the physical reason behind the body acceleration?} The answer is not univocal, for instance one could think about some perturbation of the fluid or, alternatively, to the action of some external forces. In this paper, we focus on the first case \par Another contribution in this direction is the work by Hishida and Maremonti \cite{His-Mar}, who considered the situation in which the fluid is subject to an initial perturbation in $L^3(\OO)$. Partially following the ideas introduced in \cite{Mar-L3Lor}, they proved the attainability of the steady state as the asymptotic limit of a weak solution. \par The choice of an initial perturbation in $L^3(\OO)$ is due to the fact that a {\it physically reasonable} stationary regime does not have, {\it a priori}, a finite kinetic energy, see \cite{Finn-sta, Ga:Book}. However, it is also known that a  {\it physically reasonable} solution to the stationary Navier-Stokes problem in exterior domains belongs to $L^q(\OO)$, for all $q>2$. This aspect motivates our research and the requirement for an initial perturbation of the fluid in $L^p(\OO)$, $p\in (2,3)$. \par Although, as in \cite{His-Mar}, the attainability of the stationary regime is deduced through the construction and the investigation of the time-asymptotic properties of a weak solution, the argument lines used here seem to be more immediate with respect to those in \cite{His-Mar}. Of course, they can be adapted to the situation proposed in the quoted paper. \par We essentially follow an idea employed by the author in the paper \cite{Pal}, where the existence and the partial regularity of a weak solution to the Navier-Stokes initial boundary problem was obtained, extending and improving a result due to Calder\'on \cite{Cald} (see section\,\ref{sec: SOP} for more details). \par In analogy to the paper \cite{His-Mar}, we might also inquire about the reverse situation, the {\it Landing Problem}. In this case, the body stops after some finite time and one has to investigate the asymptotic-in-time convergence of the fluid to the rest state. Since in that case we have in mind a different decomposition approach, we postpone the analysis of this problem to a future work. \par Let us conclude this brief overview by pointing out that in two dimensions the {\it Starting Problem} seems to be an intriguing open question. However, there are results related to the stability of a {\it physically reasonable} steady state, see \cite{Mae}. 

\subsection{Strategy of the proof}\label{sec: SOP}
Our result is achieved in several steps. Firstly, generalizing an idea employed by Calder\'{o}n in \cite{Cald} in the context of the Navier-Stokes Cauchy problem, and extended  to the case of the initial boundary value problem in \cite{Pal} (see also \cite{Pop} for an extension to Besov spaces), we decompose the initial velocity of the fluid as the sum
\[
u_0=v_0+w_0\,,\,\,v_0\in L^3(\OO)\,\,\text{and }w_0\in L^2(\OO)\,.
\]
The weak solution $u$ is obtained as the sum $u=v+w$, where $(v,\pi_v)$ is a regular solution to
\begin{equation}\label{aux-start}
 \begin{cases}
        v_t-\Delta  v+v\cdot\n v+hu_\infty\cdot\n v+\n\pi_v=0\,,\,\,\text{in }(0,T)\times\OO\,,\\ \hskip5.3cm\n\cdot v=0\,,\,\,\text{in }(0,T)\times\OO\,,\\ \hskip3.75cmv(t,x)=-h(t)u_\infty\,,\,\,\text{on }(0,T)\times\partial\OO\,, \\ \hskip2.8cm\displ\lim_{|x|\to\infty }v(t,x)=0\,,\\ \hskip2.8cm(v,\varphi)=(v_0,\varphi)\,,\,\,\text{for }t=0\,,\,\,\text{for all }\varphi\in\mathscr C_0(\OO)\,,
    \end{cases}
\end{equation}
while $w$ is a weak solution to
\begin{equation}\label{pert-start}
     \begin{cases}
        w_t-\Delta  w+w\cdot\n w+w\cdot\n v+v\cdot\n w+hu_\infty\cdot\n w+\n\pi_w=0\,,\\ \hskip4cm\n\cdot w=0\,,\,\,\text{in }(0,T)\times\OO\,,\\ \hskip3.77cmw(t,x)=0\,,\,\,\text{on }(0,T)\times\partial\OO\,, \\ \hskip2.84cm\displ\lim_{|x|\to\infty }w(t,x)=0\,,\\ \hskip2.8cm(w,\varphi)=(w_0,\varphi)\,,\,\,\text{for }t=0\,,\,\,\text{for all }\varphi\in\mathscr C_0(\OO)\,.
    \end{cases}
\end{equation}
The analysis of problems \eqref{aux-start}-\eqref{pert-start} is carried out by different techniques. \par As far as problem \eqref{aux-start} is concerned, the existence of a unique regular solution is obtained by suitably adapting the arguments employed by Galdi et al. in \cite{GHS} to the case of a ``small" $L^3$ initial datum. Hence, we follow an evolution equations approach and a fixed point argument. In order to employ the quoted technique, we need a divergence-free initial datum. However, in the more general case of an initial condition of the kind \eqref{aux-start}$_5$, by  virtue of the Helmholtz decomposition due to Solonnikov \cite{Sol}, see also \cite{Mar-revSol}, we find that the solution corresponding to the initial datum $\mathbb P_3 v_0$ also solves \eqref{aux-start}. \par The uniqueness of the solution to problem \eqref{aux-start} is obtained, in $[0,T_0]$, through a duality argument, firstly introduced by Foias in \cite{Foi} in the context of the Navier-Stokes equations. For $t> T_0$, problem \eqref{aux-start} reduces to the exterior Navier-Stokes problem investigated by Shibata in \cite{Shi}, hence the uniqueness holds in $(T_0,\infty)$ in view of the results therein. \par Problem \eqref{pert-start} essentially describes the perturbation of a fluid motion with respect to a nonstationary state, with an initial perturbation in $L^2$. Hence, we follow an $L^2$-theory approach, as done by Leray in the case of a perturbation from the rest state. We first construct a weak solution as the weak limit, in the metric detected by the energy relation, of a suitable approximating sequence. Then, we discuss the partial regularity of the weak solution in terms of {\it Th\'eor\`eme de Structure}. For this task, we follow the argument lines proposed by Leray in \cite{L}. However, an alternative proof based on the arguments developed in \cite{Mar-Pal-NS, Mar-Pal-FSI}, and also employed in \cite{Mar-TDS,Pal}, is certainly possible. \par Therefore, the well-posedness results related to problems \eqref{aux-start}-\eqref{pert-start} allow us to obtain a  weak solution $u$ to problem \eqref{NS-Ose}, in the sense of Definition\,\ref{def: WS-Ose}, as the sum $u=v+w$. \par In particular, as we find that our weak solution is regular for large times, we are allowed to discuss, asymptotically in time, the attractiveness of the steady-state obtained as the solution to problem \eqref{NS-Ose-staz}. 
\section{Notation and Preliminary results}\label{sec: prel}
We set
\[
\mathscr C_0(\OO):=\{\varphi\in C_0^\infty(\OO)\,:\,\n\cdot\varphi=0\}\,.
\]
We denote by $L^p(\OO)$, $W^{m,p}(\OO)$ and $\widehat W^{m,p}(\OO)$, $m\in\N$, $p\in [1,\infty]$, respectively, the Lebesgue, Sobolev and homogeneous Sobolev spaces. By the symbol $L^{p,q}(\OO)$, $p,q\in[1,\infty]$, we refer to the Lorentz spaces.\par 
We define the spaces $J^p(\OO)$ and $J^{m,p}(\OO)$ as the closures of $\mathscr C_0(\OO)$ with respect to the metrics $\dm \cdot\dm_p$ and $\dm \cdot\dm_{m,p}$, respectively. It is known that, for all $p\in (1,\infty)$, the following decomposition holds:
\[
L^p(\OO)=J^p(\OO)\oplus G^p(\OO)\,,
\]
where $G^p(\OO)$ denotes the subspace of $L^p(\OO)$ of functions of the kind $f=\n h$, $h\in W^{1,p}_{loc}(\OO)$ with $\n h\in L^p(\OO)$, see \cite{Ga:Book, Sol}. Moreover, we are able to define a bounded linear operator
\[
\mathbb P_p:L^p(\OO)\to J^p(\OO)\,,
\]
known as Helmholtz projection.\par  Given $-\infty\leq a < b\leq+\infty$, we denote by  $C(a,b;X)$ and $L^p(a,b;X)$, $X$ Banach space, the Bochner spaces. We will usually employ the concise notation $(u,\phi)$ to denote the integral 
\[
\int_\OO u\cdot\phi\,dx\,.
\]
 By the symbol $W^{m-\frac1p,p}(\partial\OO)$, we refer to the trace space of $W^{m,p}(\OO)$. If $p=2$, we might employ the alternative notation $H^{m-\frac12}(\partial\OO)$. The symbol $\gamma$ is employed for the usual Gagliardo trace operator.\par  For an overview on function spaces we refer the reader to \cite{Adams}, while for the spaces of the hydrodynamics we quote \cite{Ga:Book}. \par Following the notation employed in \cite{Ga:Book}, by $B_R$ we denote the ball centered in the origin with radius $R$, and we set $\OO_R=\OO\cap B_R$ and $\OO^R=\OO\setminus B_R$. \par In the following, by $c$ we mean a generic constant, whose value may vary from time to time and is inessential for our purposes.\par We now state some lemmas that will be needed for our goals.
\begin{lemma}\label{Dec-Lp}
    {\sl Let $p,q,r\in [1,\infty)$ be such that $r<p<q$. For all $\rho>0$ and for all $f\in L^p(\Omega)$ there exist $f_1\in L^q(\Omega)$ and $f_2\in L^r(\Omega)$ such that $f=f_1+f_2$ and 
    \be\label{dec-est}
    \ba{rl}\displ
\dm f_1\dm_q&\hskip-0.2cm\leq \rho^{1-\frac{p}{q}}\dm f\dm_p^{\frac{p}{q}}\,,\VS \dm f_2\dm_r&\hskip-0.2cm\leq \rho^{1-\frac{p}{r}}\dm f\dm_p^{\frac{p}{r}}\,.
\ea
    \ee}
\end{lemma}
\bp
The result is proved in \cite{Pal} and it is a more immediate reproduction of a decomposition result due to Calder\'on \cite{Cald}.
\ep
\begin{lemma}\label{le: Sto-interp}
    {\sl Let $g\in J^{1,2}(\OO)\cap W^{2,2}(\OO)$. Then there exists a constant $c$, independent of $g$, such that
    \[
    \dm g\dm_\infty\leq c\dm\mathbb P_2\Delta g\dm_2^\frac12\dm \n g\dm_2^\frac12\,.
    \]}
\end{lemma}
\bp
See \cite{Mar-int}.
\ep
\begin{lemma}\label{le: Mar-AB}
    {\sl Let $y(t)\in C^{1}(t_0,\infty)$ such that $y(t)\geq 0$ and assume that
    \[
    \begin{cases}
        \dot y(t)\leq a(t)y(t)+K_1y^2(t)+K_2y^3(t)\,,\\ y(t_0)=y_0\,,\\  E=\int_{t_0}^\infty y(t)\,dt<\infty\,,
    \end{cases}
    \]
    for some $K_1,K_2\ge 0$. Assume further that there exists $K_3\ge 0$ such that $a(t)\leq K_3$ for all $t\ge t_0$, and that for some $\delta>0$
    \[
    E<\frac{\delta}{2K(1+\delta)^2}\,,\,\,\text{with }K:=\max\{K_i,\, i=1,\dots,3\}\,.
    \]
    Then $$y(t_0)\leq \frac\delta2\implies y(t)\leq \delta\,,\,\,\text{for all }t\ge t_0\,.$$}
    \end{lemma}
    \bp
The proof of the lemma can be found in \cite{Mar-AB}.
    \ep
We recall a well-posedness result related to the steady-state problem \eqref{NS-Ose-staz}.
\begin{lemma}\label{le: WP-staz}
    {\sl There exists $\sigma_0>0$ such that: for all $u_\infty\in\R^3\setminus\{0\}$ with $|u_\infty|\leq \sigma_0$, problem \eqref{NS-Ose-staz} admits a unique solution $(U,\pi_U)$ such that
    \begin{equation}\label{est-staz}
    \ba{rl}
\dm U\dm_q&\hskip-0.2cm\leq c(q)|u_\infty|^\frac12\,,\,\,\,\,\text{for all }q\in (2,\infty]\,,\VS \dm \n U\dm_r&\hskip-0.2cm\leq c'(r)|u_\infty|^\frac12\,,\,\,\text{for all }r\in (\frac43,\infty]\,,\VS \dm \pi_U\dm_s&\hskip-0.2cm\leq c''(s)|u_\infty|^\frac12\,,\,\text{for all }s\in(3,\infty)\,.
    \ea
    \end{equation}
    Moreover, $U$ satisfies, as $|x|\to\infty$, 
    \begin{equation}\label{As-staz-Finn}
        |U(x)|=O(|x|^{-\frac12-\vep})\,,
    \end{equation}
    for some $\vep>0$.}
\end{lemma}
\bp
The existence and uniqueness of a solution, as long as property \eqref{As-staz-Finn}, are due to Finn \cite{Finn-sta}, see also \cite[Chapter X]{Ga:Book}. Estimates \eqref{est-staz} follow easily from \cite[Theorem 1.1]{Shi}.
\ep
Let $(U,\pi_U)$ be the solution to problem \eqref{NS-Ose-staz} ensured by Lemma\,\ref{le: WP-staz}. Let us consider consider the following initial boundary value problem
    \begin{equation} \label{pert-staz}
    \begin{cases}
       \overline u_t-\Delta \overline u+\overline u\cdot\n\overline u+u_\infty\cdot\n\overline u+U\cdot\n\overline u+\overline u\cdot\n U+\n\pi_{\overline u}=0\,,\,\,\text{in }(0,T)\times\OO\,,\\ \hskip3.99cm\n\cdot\overline u=0\,,\,\,\text{in }(0,T)\times\OO\,,\\ \hskip3.75cm\overline u(t,x)=0\,,\,\,\text{on }(0,T)\times\partial\OO\,, \\ \hskip2.8cm\displ\lim_{|x|\to\infty }\overline u(t,x)=0\,,\\ \hskip4.5cm\overline u=\overline u_0\,,\,\,\text{on }\{0\}\times\OO\,.
    \end{cases}
\end{equation}
Denoting by $T_{u_\infty}$ the Oseen semigroup related to the vector $u_\infty$, see \cite{Kob-Shi, Miy, Shi}, we convert problem \eqref{pert-staz} into the abstract integral equation
\begin{equation}\label{mild-pert-staz} 
    \overline u(t)=T_{u_\infty}(t)\overline u_0-\int_0^tT_{u_\infty}(t-\tau)\mathbb P[\overline u\cdot\n \overline u+\overline u\cdot\n U+U\cdot\n\overline u](\tau)\,d\tau\,.
\end{equation}
A function $\overline u$ satisfying the integral equation \eqref{mild-pert-staz} is called a {\it mild solution} to problem \eqref{pert-staz}. \par We state the following well-posedness result.
\begin{lemma}\label{le: WP-Shibata}
    {\sl There exists $\vep>0$ such that for all $u_\infty\in\R^3\setminus\{0\}$ and $\overline u_0\in J^3(\OO)$ satisfying $|u_\infty|+\dm \overline u_0\dm_3\leq \vep$, problem \eqref{pert-staz} admits a unique mild solution $\overline u\in C([0,\infty);J^3(\OO))$ verifying the following properties:
    \begin{equation}\label{reg-Shi}
        t^{\frac32(\frac13-\frac1q)}\overline u\in C([0,\infty);L^q(\OO)),\,\,\forall q\in [3,\infty];\,\,\,t^\frac12\n\overline u\in C([0,\infty);L^3(\OO)).
    \end{equation}}
\end{lemma}
\bp
The proof is given in \cite{Eno-Shi,Shi}.
\ep
\section{Proof of the main result: Theorem\,\ref{thm: MR}}
In this section we prove Theorem\,\ref{thm: MR}. Employing Lemma\,\ref{Dec-Lp}, we can write $u_0$ as the sum $u_0=v_0+w_0$, with $v_0\in L^3(\OO)$ and $w_0\in L^2(\OO)$. Hence, we construct the weak solution $u$, stated in Theorem\,\ref{thm: MR}, as the sum $u=v+w$, with $v$ and $w$ that solve, respectively, problem \eqref{aux-start} and  problem \eqref{pert-start}. \par 
The quoted subproblems are investigated in details in the next subsections.  
\subsection{Analysis of problem \eqref{aux-start}}\label{sec: auxprob}
The aim of this subsection is to prove the global existence and uniqueness of a regular solution to problem \eqref{aux-start}. \par For the construction of the quoted solution we follow the argument lines given in \cite{Eno-Shi, GHS, His-Mar, Shi}. Hence, we obtain $(v,\pi_v)$ as the sum 
\[
v(t,x)=\overline v(t,x)+h(t)U(x)\,,\,\,\pi_v(t,x)=\pi_{\overline v}(t,x)+h(t)\pi_U(x)
\]
where $(U,\pi_U)$ is the unique solution to problem \eqref{NS-Ose-staz}, given in Lemma\,\ref{le: WP-staz}, and, recalling Assumption\,\ref{hfunct} and setting
\[
f(t,x):=-h'(t)U(x)+(h(t)-h^2(t))(U(x)+u_\infty)\cdot\n U(x)\,,
\]
$(\overline v,\pi_{\overline v})$ solves
\begin{equation}\label{NS-Ose-diff}
    \begin{cases}
        \overline v_t-\Delta \overline v+\overline v\cdot \n \overline v+h(t)(\overline v\cdot \n U+U\cdot\n\overline v+u_\infty\cdot \n \overline v)+\n \pi_{\overline v}=f,\\\hskip6cm\n\cdot \overline v=0,\,\,\quad \text{in }(0,T)\times\OO,\\ \hskip5.78cm\overline v(t,x)=0,\,\,\quad\text{on }(0,T)\times\partial\OO,\\\hskip4.85cm \displ \lim_{|x|\to\infty}\overline v(t,x)=0,\\\hskip4cm (\overline v,\varphi)=(v_0,\varphi),\,\,\text{for }t=0\,,\,\,\text{for all }\varphi\in\mathscr C_0(\OO).
    \end{cases}
\end{equation}
Let $v_0\in J^3(\OO)$. We convert problem \eqref{NS-Ose-diff} into the integral equation
\begin{equation}\label{Int-OseNS-diff}
    \overline v(t)=T_{u_\infty}(t)v_0+H_U[\overline v]-L_U[\overline v]-N[\overline v]\,,
\end{equation}
with
\[\hskip-0.15cm
\ba{rl}
H_U[\overline v(t)]&\hskip-0.25cm\displ:=\!\int_0^t\!T_{u_\infty}(t\!-\!\tau)\mathbb P_3[f(\tau)\!+\!(1\!-\!h(\tau))u_\infty\!\cdot\!\n \overline v(\tau)\!+\!(1\!-\!h(\tau))U\!\cdot\!\n \overline v(\tau)\!+\!(1\!-\!h(\tau))\overline v(\tau)\!\cdot\!\n U]d\tau,\VS 
L_U[\overline v(t)]&\hskip-0.25cm:=\!\displ\int_0^t T_{u_\infty}(t\!-\!\tau)\mathbb P_3[U\!\cdot\!\n \overline v(\tau)+\overline v(\tau)\!\cdot\!\n U]d\tau,\VS N[\overline v(t)]&\hskip-0.25cm:=\!\displ\int_0^t T_{u_\infty}(t\!-\!\tau)\mathbb P_3[\overline v(\tau)\!\cdot\!\n \overline v(\tau)]d\tau.
\ea
\]
Our goal is to prove the following result.
\begin{tho}\label{thm: GWP-aux}
    {\sl Let $h$ satisfy Assumption\,\ref{hfunct}. There exists $\delta>0$ such that: for all $v_0\in L^3(\OO)$ and $u_\infty\in \R^3\setminus\{0\}$, with $\dm \mathbb P_3 v_0\dm_3+|u_\infty|\leq \delta$, problem \eqref{aux-start} admits a unique solution $( v,\pi_{ v})$ with $(\overline v,\pi_{\overline v}):=(v-hU,\pi_v-h\pi_U)$ satisfying
    \begin{equation}\label{reg-sol-diff}
        \ba{c}
     t^{\frac32(\frac13-\frac1q)}\overline v\in L^\infty(0,\infty;L^q(\OO)),\,\,\text{for all }3< q\leq \infty,\VS t^\frac12\n\overline v\in L^\infty(0,\infty;L^3(\OO)),\VS \lim_{t\to 0^+}(\overline v(t),\varphi)=(v_0,\varphi),\,\,\text{for all }\varphi\in\mathscr C_0(\OO).
        \ea
    \end{equation}
    In particular, we get
    \be \label{cont-0}
    \lim_{t\to 0^+}\dm \overline v(t)-\mathbb P_3v_0\dm_3=0.
    \ee
    Moreover, for all $T>0$, it holds that
    \begin{equation}\label{sol-aux-norm}
        \sup_{(0,T)}\big[\dm \overline v(t)\dm_3+t^\frac12\dm \n\overline v(t)\dm_3+t^{\frac32(\frac13-\frac1q)}\dm\overline v(t)\dm_q \big]\leq c(\dm \mathbb P_3 v_0\dm_3+|u_\infty|^\frac12),\,\,\text{for all }q\geq 3\,.
    \end{equation}}
\end{tho}
\bp[Proof of Theorem\,\ref{thm: GWP-aux}]
The proof is essentially divided into three steps. We first prove the theorem for an initial datum $v_0\in J^3(\OO)$. Then, recalling the Helmholtz decomposition on $L^q$ spaces, we discuss the case of an initial datum in $L^3(\OO)$. Finally, we discuss the uniqueness of the solution employing a duality argument.\par We obtain $(v,\pi_v)$ as $(v,\pi_v)=(\overline v+hU,\pi_{\overline v}+h\pi_U)$, where $(U,\pi_U)$ is the solution to problem \eqref{NS-Ose-staz}, achieved in Lemma\,\ref{le: WP-staz}, and $(\overline v,\pi_{\overline v})$ solves \eqref{NS-Ose-diff}. Hence, we only focus on the solvability of problem \eqref{NS-Ose-diff}.\par 
{\bf Step 1.}  Let $v_0\in J^3(\OO)$. Then, we are allowed to convert problem \eqref{NS-Ose-diff} into the abstract integral equation \eqref{Int-OseNS-diff}. Therefore, we reduce the problem to finding a global solution to the problem \eqref{Int-OseNS-diff} that has the regularity expressed in \eqref{reg-sol-diff}-\eqref{sol-aux-norm}. For this task, we simply have to adapt the proof given in \cite{GHS} in the case $v_0=0$. Since the procedure is consolidated in the literature, we omit it here. 
\par {\bf Step 2.} We now consider the case $v_0\in L^3(\OO)$. By virtue of the Helmholtz decomposition, we have
\[
v_0=\mathbb P_3 v_0\oplus \n h\,.
\]
We consider the solution $(\overline v,\pi_{\overline v})$ corresponding to the initial datum $\mathbb P_3 v_0$, obtained via the first step of the proof. Since for all $\varphi\in\mathscr C_0(\OO)$
\[
0=\lim_{t\to 0}\dm\overline v(t)-\mathbb P_3v_0\dm_3=\lim_{t\to 0}(\overline v(t)-\mathbb P_3v_0,\varphi)\,,
\]
and $(\n h,\varphi)=0$, we conclude that
\[
\lim_{t\to 0}(\overline v(t)-v_0,\varphi)=0\,.
\]
Therefore, $(\overline v,\pi_{\overline v})$ is also a solution to problem \eqref{NS-Ose-diff} for $v_0\in L^3(\OO)$. The existence of a global solution to problem \eqref{aux-start} is achieved. \par {\bf Step 3.} We suppose that $( v_1,\pi_1)$ and $ (v_2,\pi_2)$ are solutions to problem \eqref{aux-start} corresponding to the same initial datum $v_0$, with $(\overline v_1,\overline \pi_1)=(v_1-hU,\pi_1-h\pi_U)$ and $(\overline v_2,\overline \pi_2)=(v_2-hU,\pi_2-h\pi_U)$ enjoying the regularity expressed in \eqref{reg-sol-diff}-\eqref{sol-aux-norm}. We set
\[
(\widehat v,\widehat \pi):=(v_1-v_2,\pi_1-\pi_2)=(\overline v_1-\overline v_2,\overline\pi_1-\overline\pi_2)\,.
\]
We see that the pair $(\widehat v,\widehat\pi)$ is a mild solution to
\begin{equation}\label{aux-diff}
     \begin{cases}
        \widehat v_t-\Delta \widehat v+\widehat v\cdot\n v_1+v_2\cdot\n \widehat v+hu_\infty\cdot\n\widehat v+\n\widehat\pi=0\,,\,\,\text{in }(0,T)\times\OO\,,\\ \hskip3.98cm\n\cdot\widehat v=0\,,\,\,\text{in }(0,T)\times\OO\,,\\ \hskip3.75cm\widehat v(t,x)=0\,,\,\,\text{on }(0,T)\times\partial\OO\,, \\ \hskip2.8cm\displ\lim_{|x|\to\infty }v(t,x)=0\,,\\ \hskip2.95cm(\widehat v,\varphi)=0\,,\,\,\text{for }t=0\,,\,\,\text{for all }\varphi\in\mathscr C_0(\OO)\,,
    \end{cases}
\end{equation}
We are going to employ a duality argument.
We consider $(z,\pi_z)$ solution to the nonstationary Stokes problem with initial datum $z_0\in \mathscr C_0(\Omega)$. Let $t\in (0,T_0]$. We set $\widehat z(\tau ,x)= z(t-\tau,x)$, $\tau \in (0,t)$. We multiply equation \eqref{aux-diff}$_1$ by $\widehat z$ and integrate over $(0,t)\times \Omega$. We get\footnote{ We are essentially claiming that a mild solution satisfies a very weak formulation of problem \eqref{aux-diff} with a test function having the regularity of $\widehat z$. This fact is well-known in the context of the Navier-Stokes equations, see \cite{Koz-Yam}, hence we do not spend more words here.}
\be \label{dual-uniq}
(\widehat v(t),z_0)=(\widehat v(0),z(t))+\int_0^t \big[(\widehat v\cdot \n \widehat z, v_1)+(v_2\cdot\n\widehat z,\widehat v)+(hu_\infty\cdot\n \widehat z, \widehat v)\big]\,d\tau \,.
\ee
By virtue of condition \eqref{aux-diff}$_5$, we have $(\widehat v(0),z(t))=0$. Moreover, we have
\[\hskip-0.5cm
\ba{l}\displ
\Bigg|\int_0^t \big[(\widehat v\cdot \n \widehat z, v_1)+(v_2\cdot\n\widehat z,\widehat v)+(hu_\infty\cdot\n \widehat z,\widehat v)\big]\,d\tau\Bigg|\leq \int_0^t\big|(\widehat v\cdot \n \widehat z, v_1)+(v_2\cdot\n\widehat z,\widehat v)+(hu_\infty\cdot\n \widehat z,\widehat v)\big|\,d\tau\VS\hskip5.4cm\leq \int_0^t\dm\widehat v(\tau)\dm_3\dm \n\widehat z(\tau)\dm_\frac32(\dm v_1(\tau)\dm_\infty+\dm v_2(\tau)\dm_\infty+|h(\tau)||u_\infty|)\,d\tau\,.
\ea
\]
Set $$ k:=\int_0^t\tau^{-\frac12}(t-\tau)^{-\frac12}\,d\tau\,,\,\,M:=\displ\max_{[0,T_0]}|h(\tau)|\,.$$ We remark that, for all $\tau>0$,
\[
\ba{l}
\dm v_1(\tau)\dm_\infty+\dm v_2(\tau)\dm_\infty\leq \dm \overline v_1(\tau)\dm_\infty+\dm\overline v_2(\tau)\dm_\infty+2|h(\tau)|\dm U\dm_\infty\VS\hskip2cm \leq 2c\tau^{-\frac12}(\dm\mathbb P_3 v_0\dm_3+|u_\infty|^{\frac12})+2cM|u_\infty|^\frac12\leq 2c(\tau^{-\frac12}+M)(2\delta^\frac12+\delta)\,.
\ea
\]
Hence, employing the classical estimates for solutions to the Stokes IBVP in an exterior domain (see, e.g., \cite[Theorem 1.2]{Mar-Sol}), we obtain
\[
\ba{l}\displ
\int_0^t\dm\widehat v(\tau)\dm_3\dm \n\widehat z(\tau)\dm_\frac32(\dm v_1(\tau)\dm_\infty+\dm v_2(\tau)\dm_\infty+|h(\tau)||u_\infty|)\,d\tau\VS\hskip1cm\leq 3c(k+MT_0^\frac12)\dm z_0\dm_\frac32(2\delta^\frac12+\delta)\sup_{(0,t)}\dm \widehat v(\tau)\dm_3\,.
\ea
\]
Since there is no lower bound for $\delta$, we may assume that
\[
\overline c:=3c(k+MT_0^\frac12)(2\delta^\frac12+\delta)<1\,,
\]
Going back to \eqref{dual-uniq}, in view of the arbitrariness of $z_0$, we obtain
\[
(1-\overline c)\sup_{(0,t)}\dm \widehat v(\tau)\dm_3\leq 0\,,
\]
which in turn implies that $\widehat v\equiv 0$ on $[0,T_0]$.\par For $t> T_0$, system \eqref{aux-start} reduces to the exterior Navier-Stokes problem investigated in \cite{Shi}, hence the uniqueness for $t> T_0$ follows from the results in the quoted paper. The theorem is completely proved.
\ep
\subsection{Analysis of problem \eqref{pert-start}}\label{sec: pertprob}
We now focus on problem \eqref{pert-start}. 
\begin{rem}\label{rem: L2-costanti}
    {\rm Recalling Lemma\,\ref{Dec-Lp}, as in Theorem\,\ref{thm: GWP-aux} there is no upper bound on $\delta$, we can take it so that $A:=c_1c(\dm\mathbb P_3 v_0\dm_3+|u_\infty|^\frac12)<1$, where $c(\dm\mathbb P_3 v_0\dm_3+|u_\infty|^\frac12)$ is the right-hand side in relation \eqref{sol-aux-norm} and $c_1$ is the constant related to the Sobolev embedding $\widehat W^{1,2}(\OO)\hookrightarrow L^6(\OO)$.}
\end{rem}
We give the notion of weak solution to problem \eqref{pert-start}.
\begin{defi}\label{def: WS-pert-Ose}
    {\sl Let $w_0\in L^2(\OO)$. A field $w:(0,\infty)\times\Omega\to\R^3$ is said a weak solution to problem \eqref{pert-start} if, for all $T>0$, the following properties hold:
    \begin{description}
        \item[(a)] $w\in L^\infty(0,T;L^2(\Omega))\cap L^2(0,T;J^{1,2}(\Omega))$;
        \item[(b)] for all $t,s\in (0,T)$, $w$ satisfies the following integral equation:
        \[
        \ba{l}\displ
        \int_s^t \big[(w,\phi_\tau)-(\n w,\n\phi)+((w+v+hu_\infty)\cdot\n\phi, w)+(w\cdot\n\phi,v)\big]\,d\tau\VS\hskip0.5cm+ (w(s),\phi(s))=(w(t),\phi(t))\,,\,\,\text{for all }\phi\in C^1_0([0,T)\times \Omega),\,\text{with }\n\cdot\phi=0,\,\text{for all }t\in [0,T);
        \ea
        \]
        \item[(c)] $w$ satisfies the following inequality:
        \be\label{SEI-start}
        \dm w(t)\dm_2^2+2(1-A)\int_s^t \dm \n w(\tau)\dm_2^2\,d\tau\leq \dm w(s)\dm_2^2\,,\,\,\text{for all }t\ge s\,\,\text{and for a. a. }s\in( 0,t)\,;
        \ee
        \item[(d)] $\displ\lim_{t\to 0^+}(w(t),\varphi)=(w_0,\varphi)$, for all $\varphi\in\mathscr C_0(\Omega)$.
    \end{description}}
\end{defi}
We also give the following definition of regular solution.
\begin{defi}\label{def: RSOL-pert}
    {\sl A weak solution $w$ to problem \eqref{pert-start} is said regular on $(0,T)\times\OO$ if 
    \[
    D^2 w,\,w_t\in L^2(0,T;L^2(\Omega))
    \]
    and there exists a scalar field $\pi_w:(0,T)\times\OO\to\R$, with
    \[
    \n\pi_w\in L^2(0,T;L^2(\OO))\,,    \] such that $(w,\pi_w)$ solves \eqref{pert-start} a. e. in $(0,T)\times\OO$.}
\end{defi}
Our task is to prove the following theorem.
\begin{tho}\label{thm: L2-WPpert}
    {\sl For all $w_0\in L^2(\OO)$ there exists a weak solution $w$ to problem \eqref{pert-start}. Moreover, for all $\eta>0$ there exist $\theta\in[T_0,T_0+\eta^{-2}(1-A)^{-1}\dm \mathbb P_2w_0\dm_2^4]$, and a sequence of open time intervals $\{(\theta_\ell,\Theta_\ell)\}_{\ell\in\N_0}$, with $(\theta_0,\Theta_0)=(\theta,\infty)$, such that:
    \[
    \ba{c}
    (0,\infty)\setminus \underset{\ell\in\N_0}\cup (\theta_\ell,\Theta_\ell)\,\,\text{has zero Lebesgue measure, }\VS w\,\text{is a regular solution on }(\theta_\ell,\Theta_\ell)\times\OO,\,\text{for all }\ell\in\N_0\,.
    \ea
    \]}
\end{tho}
We achieve the proof of Theorem\,\ref{thm: L2-WPpert} in two steps: we first prove the existence of a weak solution and then discuss its partial regularity. \par The existence of a weak solution is obtained through an approximation procedure, following the ideas introduced by Leray in \cite{L}. Hence, we first study the following approximating problem. We set $\mathbb J_n:=J_\frac1n$, where $J_\frac1n$ is a Friedrich mollifier in the spatial variables.
\begin{equation}\label{Moll-pb}
    \begin{cases}
        w^n_t-\Delta w^n+\mathbb J_n[w^n]\cdot\n w^n+w^n\cdot\n v+(v+hu_\infty)\cdot\n w^n+\n\pi_{w^n}=0\,,\\
   \hskip3.13cm \n\cdot w^n=0\,,\quad\,\,\,\,\,\,\quad\,\,\,\,\,\,\text{in }(0,T)\times\Omega\,,\\
   \hskip2.89cm w^n(t,x)=0\,,\quad\quad \,\,\,\,\,\,\quad\text{on }(0,T)\times\partial\Omega\,,\\\displ\hskip2cm\lim_{|x|\to\infty}w^n(t,x)=0\,,\quad\qquad\quad\text{in }(0,T)\,, \\
   \hskip3.69cm w^n=\mathbb P_2w_0\,,\,\,\qquad\text{on }\{0\}\times\Omega\,.
    \end{cases}
\end{equation}
The existence and uniqueness of a solution to problem \eqref{Moll-pb} can be reached employing consolidated methods such as, for instance, the Faedo-Galerkin method, introduced by Prodi in \cite{Pro} for the Navier-Stokes problem in a bounded domain, and the invading domains technique, proposed by Heywood in \cite{Hey}. \par Therefore, following Prodi \cite{Pro}, on each bounded subdomain $\OO_R$ an orthonormal basis of $J^2(\OO_R)$ is obtained by investigating the eigenvalue problem for the Stokes operator. The aforementioned basis allows us to introduce a Galerkin approximating sequence on each bounded subdomain $\OO_R$ and, performing some a priori estimates on these sequences that do not depend on the radius $R$, we fulfill our goal in the way suggested by Heywood in \cite{Hey}. \par In the following lemma, we collect some basic estimates that are obtained for the Galerkin approximating sequence $\{w^n_{j,R}\}_{j\in\N}$ related to the subdomain $\OO_R$. For their proof, it is sufficient to simply adapt to the present situation the arguments developed in \cite[Lemma 11, Lemma 12]{Pal}. Therefore, we do not prove them here.\par  Set $M:=\displ\max_{[0,T_0]}|h(t)|$.
\begin{lemma}\label{le: ENES-start}
    {\sl For all $n\in\N$ and $t>0$, we get, uniformly in $j\in\N$,
    \begin{equation}\label{En-Gal}
        \dm w^n_{j,R}(t)\dm_2^2+2(1-A)\int_0^t\dm \n w^n_{j,R}(\tau)\dm_2^2\,d\tau\leq \dm\mathbb P_2 w_0\dm_2^2\,, 
        \end{equation}
        and
        \begin{equation}
        \ba{l}\displ
         t\dm \n w^n_{j,R}(t)\dm_2^2 + \int_0^t \tau \big[\dm \mathbb P_2\Delta w^n_{j,R}(\tau)\dm_2^2+\dm \partial_\tau w^n_{j,R}(\tau)\dm_2^2\big]\,d\tau\leq \int_0^t\dm  \n w^n_{j,R}(\tau)\dm_2^2\,d\tau\VS\hskip0.5cm+c\int_0^t \tau\big[c(n)\dm w^n_{j,R}(\tau)\dm_2^2\dm\n w^n_{j,R}(\tau)\dm_2^2\!+\!(\dm \n v(\tau)\dm_3^2\!+\!\dm v(\tau)\dm_\infty^2\!+\!|u_\infty|^2)\dm\n w^n_{j,R}(\tau)\dm_2^2 \big]\,d\tau\,,
         \ea
    \end{equation}
    for some $c=c(M)$.}
\end{lemma}
Lemma\,\ref{le: ENES-start} allows us to state the following existence and uniqueness result for problem \eqref{Moll-pb}. We do not give the proof as it is based on a classical limit process and on a straightforward application of the Gronwall Lemma.
\begin{lemma}\label{le: WP-moll-start}
    {\sl For all $w_0\in L^2(\OO)$ problem \eqref{Moll-pb} admits a unique solution $(w^n,\pi_{w^n})$ such that, for all $T>0$,
    \[
    \ba{c}\displ
    w^n\in L^\infty(0,T;J^2(\OO))\cap L^2(0,T;J^{1,2}(\OO))\,, \VS \text{for all }\delta\in (0,T)\,,\, w^n\in L^\infty(\delta,T;W^{1,2}(\OO))\cap L^2(\delta,T;W^{2,2}(\OO))\,,\, \n\pi_{w^n},w^n_t\in L^2(\delta,T;L^2(\OO))\,.
    \ea
    \]}
\end{lemma}
In order to prove Theorem\,\ref{thm: L2-WPpert}, we need to prove for the approximating sequence $\{w^n\}_n$, ensured by Lemma\,\ref{le: WP-moll-start}, the following convergence property:
\[
\text{for all }t>0\,,\,\,\lim_{R\to\infty}\dm w^n(t)\dm_{L^2(\OO^R)}=0\,,\,\,\text{uniformly in }n\in\N\,.
\]
To achieve the above limit property, we need to perform some estimates on the pressure terms $\{\pi_{w^n}\}_n$, given in the following lemma.
\begin{lemma}\label{le: Press-moll}
    {\sl Let $(w^n,\pi_{w^n})$ be the solution to problem \eqref{Moll-pb} established in Lemma\,\ref{le: WP-moll-start} and let $\lambda\in (0,\frac12)$. Set $\beta=\frac{1-2\lambda}{3}$. Then, $\pi_{w^n}=\pi^n_1+\pi_2^n+\pi_3^n$, with
    \be \label{press-est-moll}
    \ba{c}\displ
    \int_0^T\dm \pi_{1}^n(\tau)\dm_\frac{15}{7}^\frac54\,d\tau\leq C_1(\dm v_0\dm_3,T)\,,\VS\dm\pi_2^n(t)\dm_{L^r(\OO^{R})}\leq C_2(T)\dm\mathbb P_2 w_0\dm_2[t^{-1+\frac\beta2}+t^{-1+\frac\beta4}]\,,\,\,\text{for all }R>0\,, \,\,r>3\,\, \,\,\text{and }t\in (0,T)\,,\VS \int_0^T \dm\pi_3^n(\tau)\dm_6^2\leq C_3(|u_\infty|, \dm w_0\dm_2,T)\,. 
    \ea
    \ee
    }
    \end{lemma}
    \bp
    We can write $(w^n,\pi_{w^n})=(w^n_1,\pi_1^n)+(w_2^n,\pi^n_2)+(w^n_3,\pi^n_3)$, where $(w^n_1,\pi^n_1)$ solves 
    \begin{equation}\label{Sto-w1}
        \begin{cases}
            \partial_tw^n_1-\Delta w^n_1=-\n\pi_1^n+f^n\,,\quad \text{in }(0,T)\times\OO\,,\\
            \hskip1cm\n\cdot w^n_1=0\,,\quad\qquad\quad\,\,\,\quad\text{in }(0,T)\times\OO\,,\\
            \hskip0.8cmw^n_1(t,x)=0\,,\quad\qquad\quad\,\,\,\quad\text{in }(0,T)\times\partial\OO\,,\\ \hskip1.63cmw^n_1=0\,,\quad\qquad\quad\,\,\,\quad\text{on }\{0\}\times\OO\,\,,
        \end{cases}
    \end{equation}
    with
    \[
    f^n=-\mathbb J_n[w^n]\cdot\n w^n-v\cdot\n w^n-w^n\cdot\n v\,,
    \]
 $(w^n_2,\pi^n_2)$ solves
    \begin{equation}\label{Sto-w2}
        \begin{cases}
            \partial_tw^n_2-\Delta w^n_2=-\n\pi_2^n\,,\,\,\quad\quad\quad \text{in }(0,T)\times\OO\,,\\
            \hskip1cm\n\cdot w^n_2=0\,,\quad\qquad\quad\,\,\,\quad\text{in }(0,T)\times\OO\,,\\
            \hskip0.8cmw^n_2(t,x)=0\,,\quad\qquad\quad\,\,\,\quad\text{in }(0,T)\times\partial\OO\,,\\ \hskip1.63cmw^n_2=\mathbb P_2 w_0\,,\quad\,\,\quad\,\,\,\quad\text{on }\{0\}\times\OO\,\,,
        \end{cases}
    \end{equation}
    and $(w^n_3,\pi^n_3)$ solves
    \begin{equation}\label{Sto-w2}
        \begin{cases}
            \partial_tw^n_3-\Delta w^n_3=-\n\pi_3^n+g\,,\,\,\quad\,\,\text{in }(0,T)\times\OO\,,\\
            \hskip1cm\n\cdot w^n_3=0\,,\quad\qquad\quad\,\,\,\quad\text{in }(0,T)\times\OO\,,\\
            \hskip0.8cmw^n_3(t,x)=0\,,\quad\qquad\quad\,\,\,\quad\text{in }(0,T)\times\partial\OO\,,\\ \hskip1.67cmw^n_3=0\,,\quad\,\,\quad\,\,\,\quad\quad\,\,\,\text{on }\{0\}\times\OO\,\,,
        \end{cases}
    \end{equation}
    with 
    \[
    g:=-hu_\infty\cdot\n w^n\,.
    \]
    Estimates \eqref{press-est-moll}$_{1,2}$ are proved like in \cite[Lemma 13, Remark 7]{Pal}. We prove estimate \eqref{press-est-moll}$_3$. Actually, since $g\in L^2(0,T;L^2(\OO))$, we get
    \[
    \int_0^T\dm \n \pi^n_3(\tau)\dm_2^2\,d\tau\leq C(T)\int_0^T\dm g(\tau)\dm_2^2\,d\tau \leq C(T,|u_\infty|,\dm\mathbb P_2w_0\dm_2,\max_{[0,T]}|h(\tau)|)\,.
    \]
    Hence, by re-defining $\pi_3^n$ up to a suitable function of time, we conclude the proof.
    \ep
    \begin{rem}
        {\rm We point out that we also have $\n\pi_2^n\in L^2(0,T;L^2(\OO))$. This regularity property for $\pi^n_2$ would be sufficient in order to prove the next Lemma. However, we found it worth to state some stronger regularity properties for this contribution to the pressure field.}
    \end{rem}
We can now state the following convergence result.
\begin{lemma}\label{le: Conv-extball}
        {\sl Let $(w^n,\pi_{w^n})$ be the solution to problem \eqref{Moll-pb} established in Lemma\,\ref{le: WP-moll-start}. Then, for all $t>0$, we have
        \be
        \ba{c}
\dm w^n(t)\dm_{L^2(\OO^{2R})}\leq \dm \mathbb P_2 w_0\dm_{L^2(\OO^R)}+c(t,\dm v_0\dm_3,\dm \mathbb P_2 w_0\dm_2)\Psi(R)\,,\,\, \text{and }\displ\lim_{R\to\infty}\Psi(R)=0\,.
\ea
        \ee
        In particular, we get, for all $t>0$,
        \begin{equation}
            \lim_{R\to\infty}\dm w^n(t)\dm_{L^2(\OO^R)}=0\,,\,\,\text{uniformly in }n\in\N\,.
        \end{equation}}
    \end{lemma}
Before proving the lemma, we introduce a suitable smooth cut-off function. Let us consider the function
\[
\kappa(\epsilon):=\begin{cases}
    e^{-\frac{1}{\epsilon}},\,\,\text{if }\epsilon>0,\\ 0,\quad\,\,\,\text{if }\epsilon\leq 0,
\end{cases}
\]
and set 
\[
\vartheta(\varrho):=\frac{\kappa(1-\varrho)}{\kappa(\varrho)+\kappa(1-\varrho)}\,,\,\,\varrho\in\R.
\]
It is readily seen that $\vartheta\in C^\infty(\R)$ and
\[
\begin{cases}
    \vartheta(\varrho)\equiv 0,\,\,\text{if }\varrho\geq 1,\\ \vartheta(\varrho)\equiv 1,\,\,\text{if }\varrho\leq 0.
\end{cases}
\]
Let $R>1$. We consider the function $\chi_R:\R^3\to \R$ defined by
\[
\chi_R(x):=\begin{cases}
    1,\qquad\qquad\quad \,\text{if }|x|\leq R,\\ \vartheta\big(\frac{\log\big(\frac{|x|}{R}\big)}{\log R}\big),\quad\text{if }R<|x|<R^2,\\ 0,\qquad\qquad\quad \,\text{if }|x|\geq R^2.
\end{cases}
\]
We see that $\chi_R\in C_0^\infty(\R^3)$ and 
\be \label{gradcutoff}
\n \chi_R(x)=\begin{cases}
    0,\qquad\qquad\quad\qquad\quad\,\,\,\, \,\text{if }|x|\leq R,\\ \vartheta'\big(\frac{\log\big(\frac{|x|}{R}\big)}{\log R}\big)\frac{x}{|x|^2\log R},\quad\text{if }R<|x|<R^2,\\ 0,\qquad\qquad\quad \qquad\quad\,\,\,\,\text{if }|x|\geq R^2.
\end{cases}
\ee
Moreover, we also deduce
\be\label{estcutoff}
\begin{cases}
|\n\chi_R(x)|\leq c(R\log R)^{-1},\,\,\text{for all }x\in\R^3\,,\\
|D^2\chi_R(x)|\leq cR^{-2}(\log R)^{-1}\,,\,\,\text{for all }x\in\R^3\,.
\end{cases}
\ee
We now prove Lemma\,\ref{le: Conv-extball}.
\bp[Proof of Lemma\,\ref{le: Conv-extball}]
Let $R>\max\{1,\text{diam}(\R^3\setminus\OO)\}$ and set $\zeta_R(x):=1-\chi_R(x)$. We multiply equation \eqref{Moll-pb}$_1$ by $\zeta_R^2w^n$. Integration by parts over $(0,t)\times\OO$ yields
\be\label{w-extsphere}
\ba{l}\displ
\frac12\dm \zeta_Rw^n(t)\dm_2^2\!+\!\int_0^t\dm \zeta_R^2\n w^n(\tau)\dm_2^2\,d\tau=\frac12\dm \zeta_R\mathbb P_2w_0\dm_2^2+\!\int_0^t (|w^n(\tau)|^2,\Delta\zeta_R^2)\,d\tau\VS\hskip0.8cm+\!\int_0^t ([\mathbb J_n[w^n]+v(\tau)+hu_\infty]\cdot \n \zeta_R^2,|w^n(\tau)|^2)\,d\tau+\int_0^t(w^n(\tau)\cdot\n v(\tau),\zeta_R^2w^n(\tau))\,d\tau\VS\hskip4.5cm+\int_0^t(\pi_{w^n}(\tau)w^n(\tau),\n\zeta_R^2)\,d\tau=:\frac12\dm \zeta_R\mathbb P_2w_0\dm_2^2+\sum_{i=1}^4I_i(t,n,R)\,.
\ea
\ee
By virtue of Lemma\,\ref{le: ENES-start} and estimate \eqref{estcutoff}$_2$, we get
\[
|I_1(t,n,R)|\leq c(t,\dm \mathbb P_2 w_0\dm_2^2)(R(\log R)^\frac12)^{-2}\,.
\]
Concerning $I_2$, employing H\"{o}lder's inequality and estimate \eqref{estcutoff}$_1$ we have
\[
|I_2(t,n,R)|\leq c(R\log R)^{-1}\int_0^t\big[\dm w^n(\tau)\dm_3^3+(\dm v\dm_\infty+|u_\infty|)\dm w^n(\tau)\dm_2^2\big]\,d\tau\,.
\]
Since $\dm w^n\dm_3\leq \dm w^n\dm_2^\frac12\dm\n w^n\dm_2^\frac12$ and, by virtue of Theorem\,\ref{thm: GWP-aux}, $v\in L^1(0,T;L^\infty(\OO))$, invoking Lemma\,\ref{le: ENES-start}, we deduce
\[
|I_2(t,n,R)|\leq c(t,\dm\mathbb P_2w_0\dm_2,\dm v_0\dm_3 )(R\log R)^{-1}\,.
\]
Concerning $I_3$, integrating by parts we obtain
\[
I_3(t,n,R)=-\int_0^t\big[(w^n(\tau)\cdot\n w^n(\tau),\zeta_R^2v(\tau))+(w^n(\tau)\cdot (\n\zeta_R^2\otimes w^n(\tau)),v(\tau))\big]\,d\tau=:-(I_{3,1}+I_{3,2})\,.
\]
Employing H\"{o}lder's inequality and recalling Lemma\,\ref{rem: L2-costanti}, we get
\[
|I_{3,1}(t,n,R)|\leq \int_0^t \dm w^n(\tau)\dm_6\dm \zeta_R^2\n w^n(\tau)\dm_2\dm v(\tau)\dm_3\,d\tau\leq A\int_0^t\dm \zeta_R^2\n w^n(\tau)\dm_2^2\,d\tau\,.
\]
Moreover, concerning $I_{3,2}$, employing H\"{o}lder's inequality and estimate \eqref{estcutoff}$_1$, we find that
\[
|I_{3,2}(t,n,R)|\leq c(R\log R)^{-1}\int_0^t\dm w^n(\tau)\dm_6\dm\n w^n(\tau)\dm_2\dm v(\tau)\dm_3\,d\tau.
\]
Invoking again Lemma\,\ref{le: ENES-start} and Theorem\,\ref{thm: GWP-aux}, we conclude that
\[
\lim_{R\to \infty}|I_{3,2}(t,n,R)|=0\,,\quad\text{uniformly in }n\in\N\,.
\]
Finally, we discuss the term $I_4$. We have
\[
I_4(t,n,R)\!=\!\int_0^t (\pi_{w^n}(\tau)\n\zeta_R^2,w^n(\tau))d\tau=\int_0^t ((\pi_{1}^n(\tau)\!+\!\pi^n_2(\tau)\!+\!\pi^n_3(t))\n\zeta_R^2,w^n(\tau))d\tau=:I_{4,1}\!+\!I_{4,2}\!+\!I_{4,3}\,,
\]
where
\[
\ba{rl}
I_{4,1}(t,n,R)&\hskip-0.2cm:=\displ\int_0^t (\pi_{1}^n(\tau)\n\zeta_R^2,w^n(\tau))\,d\tau\,,\VS
I_{4,2}(t,n,R)&\hskip-0.2cm:=\displ\int_0^t (\pi_{2}^n(\tau)\n\zeta_R^2,w^n(\tau))\,d\tau\,, \VS
I_{4,3}(t,n,R)&\hskip-0.2cm:=\displ\int_0^t (\pi_3^n(\tau)\n\zeta_R^2,w^n(\tau))\,d\tau\,.
\ea
\]
By virtue of H\"{o}lder's inequality and Lemma\,\ref{le: ENES-start}, and recalling \eqref{gradcutoff}, we get
\[
\ba{l}\displ
\Bigg|\int_0^t (\pi_{1}^n(\tau)\n\zeta_R^2,w^n(\tau))\,d\tau\Bigg|\displ\leq \int_0^t\dm \pi^n_1(\tau)\dm_\frac{15}{7}\dm \n \zeta_R\dm_{30}\dm w^n(\tau)\dm_2\,d\tau\VS\hskip0.7cm\leq cR^{-\frac{9}{10}}(\log R)^{-1}\int_0^t\dm \pi^n_1(\tau)\dm_\frac{15}{7}\dm w^n(\tau)\dm_2\,d\tau\leq c(\dm\mathbb P_2w_0\dm_2)R^{-\frac{9}{10}}(\log R)^{-1}\int_0^t\dm \pi^n_1(\tau)\dm_\frac{15}{7}\,d\tau\,.
\ea
\]
Recalling Lemma\,\ref{le: Press-moll}, we conclude that
\[
\lim_{R\to \infty}|I_{4,1}(t,n,R)|=0\,,\quad\text{uniformly in }n\in\N\,.
\]
We consider $I_{4,2}$. Let $r\in (3,6)$. Since $\text{supp}(\n\zeta_R)\subset \{x\,:\,R<|x|<R^2\}\subset \OO^R$, employing H\"{o}lder's inequality, Lemma\,\ref{le: Press-moll} and \eqref{gradcutoff}, we get
\[
|I_{4,2}(t,n,R)|\leq \int_0^t\dm w^n(\tau)\dm_2\dm \n \zeta_R\dm_{\frac{2r}{r-2}}\dm \pi^n_2(\tau)\dm_{L^r(\OO^R)}d\tau\leq c\dm\mathbb P_2w_0\dm_2^2R^{\frac{r-6}{2r}}(\log R)^{-1}\!\int_0^t[\tau^{-1+\frac\beta2}\!+\!\tau^{-1+\frac\beta4}]d\tau\,.
\]
We recall that from Lemma\,\ref{le: Press-moll} we have $\beta>0$. Since 
\[
\int_0^t[\tau^{-1+\frac\beta2}+\tau^{-1+\frac\beta4}]\,d\tau=\frac2\beta t^{\frac\beta2}+\frac4\beta t^{\frac\beta4}\,,
\]
and $r\in (3,6)$,
we conclude that
\[
\lim_{R\to \infty}|I_{4,2}(t,n,R)|=0\,,\quad\text{uniformly in }n\in\N\,.
\]
Finally, we estimate $I_{4,3}$. In view of \eqref{gradcutoff}, we have
\[
|I_{4,3}(t,n,R)|\leq\int_0^t \dm \pi^n_3(\tau)\dm_6\dm \n \zeta_R\dm_3\dm w^n(\tau)\dm_2\,d\tau\leq c(\log R)^{-\frac23}\int_0^t\dm \pi^n_3(\tau)\dm_6\dm w^n(\tau)\dm_2\,d\tau\,.
\]
Therefore, we deduce
\[
\lim_{R\to \infty}|I_{4,3}(t,n,R)|=0\,,\quad\text{uniformly in }n\in\N\,.
\]
Going back to \eqref{w-extsphere}, we deduce
\[
\ba{l}\displ
\frac12\dm \zeta_Rw^n(t)\dm_2^2\!+(1-A)\!\int_0^t\dm \zeta_R^2\n w^n(\tau)\dm_2^2\,d\tau\leq \frac12\dm \zeta_R\mathbb P_2w_0\dm_2^2\VS\hskip7cm +c \big(\sum_{i=1,2,4}|I_i(t,n,R)|+|I_{3,2}(t,n,R)|\big)\,.
\ea
\]
As we proved that the integral terms on the right-hand side of the above inequality tend to $0$, uniformly in $n\in\N$, as $R\to\infty$, we
conclude the proof of the Lemma.
\ep

We now prove the first part of Theorem\,\ref{thm: L2-WPpert}, hence we deduce the existence of a weak solution. 
\begin{lemma}\label{le: Ex-pert-start}
    {\sl For all $w_0\in L^2(\OO)$ there exists a weak solution $w$ to problem \eqref{pert-start} in the sense of Definition\,\ref{def: WS-pert-Ose}. Moreover, $w\in C_w(0,T;L^2(\OO))$ and satisfies, for all $t>0$,
    \be\label{eniq-0}
    \dm w(t)\dm_2^2+2(1-A)\int_0^t\dm \n w(\tau)\dm_2^2\,d\tau\leq \dm\mathbb P_2 w_0\dm_2^2\,.
    \ee}
\end{lemma} 
\bp
We recover a weak solution $w$ via a limit procedure, in the metric detected by the energy relation, on the sequence $\{w^n\}_{n\in\N}$, established in Lemma\,\ref{le: WP-moll-start}.
        By virtue of estimate \eqref{En-Gal}, we deduce that the sequence $\{w^n\}_n$ is bounded in $L^2(0,T;J^{1,2}(\OO))$. Hence, there exists a subsequence, still denoted by $\{w^n\}$, which admits a weak limit $w$ in this space.\par 
        Employing consolidated arguments, see \cite{Ga: LecNot, Lad, Pro}, one proves that, for all $\psi\in J^2(\OO)$, the sequence
        \[
        \{(w^n(\tau),\psi)\}_{n\in\N}
        \]
        is equicontinuous and uniformly bounded. Hence, the Ascoli-Arzel\`a Theorem ensures that there exists a subsequence converging to a function $k(t)\in C([0,T))$. Moreover, using well-known procedures, one deduces
        \[
        k(t)=(w(t),\psi)\,,\,\, \forall t\in[0,T)\,.
        \]
         We now show that $w$ satisfies the strong energy inequality given in the item $\bf(c)$ of Definition\,\ref{def: WS-pert-Ose}. 
        It is readily seen that the sequence $\{w^n\}_n$ satisfies the energy relation
        \be \label{SEI-moll}
        \ba{c}\displ
        \dm w^n(t)\dm_2^2+2(1-A)\int_s^t \dm\n w^n(\tau)\dm_2^2\,d\tau\leq \dm w^n(s)\dm_2^2\,, \VS\text{for all }t>s,\,\,\text{a. a. }s\in (0,t)\,\,\text{and }s=0.
        \ea
        \ee
        As a matter of fact, the above inequality is obtained multiplying equation \eqref{Moll-pb}$_1$ by $w^n$, integrating over $(s,t)\times\OO$, employing H\"{o}lder's inequality, and recalling Remark\,\ref{rem: L2-costanti}.\par 
        We show that for a. a. $s\in (0,T)$
        \[
        \lim_{n\to\infty} \dm w^n(s)\dm_2^2=\dm w(s)\dm_2^2\,.
        \]
       Recalling that, by virtue of Friedrich's Lemma and Lemma\,\ref{le: Conv-extball} we have $$\ba{c} \displ w^n\to w\,\,\text{strongly in }L^2(0,T;L^2(\OO_R))\,, \VS\lim_{R\to\infty}\dm w^n(t)\dm_{L^2(\OO^{R})}=0\,, \,\,\text{uniformly in }n\in\N\,,\,\,\text{for all }t\in(0,T)\,,
        \ea$$ 
we get
\[
\{w^n\}\to w\,\,\text{strongly in }L^2(0,T;L^2(\OO))\,,
\]
and, in particular, $\displ\lim_{n\to\infty} \dm w^n(s)\dm_2^2=\dm w(s)\dm_2^2$ a. e. in $(0,T)$.\par
        We take the minimum limit in \eqref{SEI-moll} and we get
        \[
        \dm w(t)\dm_2^2+2(1-A)\int_s^t\dm \n w(\tau)\dm_2\,d\tau\leq \dm w(s)\dm_2^2\,,\,\,\text{for all }t>0\,\text{and a. a. }s\in (0,t)\,,
       \]
       and we also deduce \eqref{eniq-0}.
        We conclude that $w$ satisfies item $\bf (c)$ of Definition\,\ref{def: WS-pert-Ose}. From the validity of item $\bf (c)$ and \eqref{eniq-0}, we also see that item $\bf (a)$ is satisfied. The proof of items $\bf (b)$ and $\bf (d)$ can be achieved via the arguments employed, for instance, in \cite{Lad}. 
        \ep
In order to prove the structure theorem for our weak solution, we need a preliminary result. \par For all $t>0$, we denote by $\mathcal T_t$ the set of instants $s\in [0,t)$ such that inequality \eqref{SEI-start} holds for the weak solution $w$ established in Lemma\,\ref{le: Ex-pert-start}. \par Set $\widetilde{\mathcal T}_t:=\mathcal T_t\cap\{s\in (0,t)\,:\, w(s)\in J^{1,2}(\OO)\}$. It is clear that $(0,t)\setminus \widetilde{\mathcal T}_t$ has zero Lebesgue measure.
\begin{lemma}\label{le: Uni-LS}
    {\sl Let $w$ be the solution established in Lemma\,\ref{le: Ex-pert-start}. Let $t>0$ and $s\in \widetilde{\mathcal T}_t$. Then, there exists $T(s)> 0$ such that 
    \[
    w\in C(s,s+T(s);J^{1,2}(\OO))\cap L^2(s,s+T(s);W^{2,2}(\OO))\,,\,\,w_t\in L^2(s,s+T(s);L^2(\OO))\,.
    \]}
\end{lemma}
\bp
Let us consider the approximating problem
\begin{equation}\label{Moll-pb-uni}
    \begin{cases}
      \overline  w^n_t-\Delta\overline w^n+\mathbb J_n[\overline w^n]\cdot\n\overline w^n+\overline w^n\cdot\n v+(v+u_\infty)\cdot\n\overline w^n+\n\pi_{\overline w^n}=0\,,\\
   \hskip3.13cm \n\cdot\overline w^n=0\,,\quad\,\,\,\,\,\,\quad\,\,\,\,\,\,\text{in }(s,T)\times\Omega\,,\\
   \hskip2.89cm\overline w^n(t,x)=0\,,\quad\quad \,\,\,\,\,\,\quad\text{on }(s,T)\times\partial\Omega\,,\\
   \hskip3.74cm\overline w^n=w(s)\,,\,\,\qquad\,\text{on }\{s\}\times\Omega\,.
    \end{cases}
\end{equation}
By following analogous procedures to those leading to the proof of Lemma\,\ref{le: WP-moll-start}, we find that problem \eqref{Moll-pb-uni} has a unique solution $(\overline w_n,\pi_{\overline w^n})$ such that, for all $T>0$,
\[
\overline w^n\in C(s,T;J^{1,2}(\OO))\cap L^2(s,T;W^{2,2}(\OO)),\, \overline w^n_t,\n\pi_{\overline w^n}\in L^2(s,T;L^2(\OO))\,.
\]
We rewrite equation \eqref{Moll-pb-uni}$_1$ as
\[
\overline w^n_t-\mathbb P_2\Delta \overline w^n=\mathbb P_2[-\mathbb J_n[\overline w^n]\cdot\n \overline w^n-\overline w^n\cdot\n v -(v+u_\infty)\cdot\n \overline w^n]\,,
\]
and we evaluate the $L^2$-norm of both sides of the above identity. We have, for all $t>s$,
\[
\ba{l}
\frac{d}{dt}\dm \n\overline w^n(t)\dm_2^2+\dm\mathbb P_2\Delta \overline w^n(t)\dm_2^2+\dm\overline w^n_t(t)\dm_2^2\leq \dm \mathbb P_2[\mathbb J_n[\overline w^n]\cdot\n \overline w^n](t)\dm_2^2\VS \hskip4.5cm+\dm\mathbb P_2[\overline w^n\cdot\n v](t)\dm_2^2+\dm\mathbb P_2[v\cdot\n \overline w^n](t)\dm_2^2+\dm\mathbb P_2[u_\infty\cdot\n \overline w^n](t)\dm_2^2\,.
\ea
\]
Recalling Lemma\,\ref{le: Sto-interp} and Theorem\,\ref{thm: GWP-aux}, employing the elementary properties of mollifiers, H\"{o}lder and Sobolev inequalities, we deduce
\[
\ba{rl}
\dm \mathbb P_2[\mathbb J_n[\overline w^n]\cdot\n \overline w^n](t)\dm_2^2&\hskip-0.2cm\leq \dm \overline w^n(t)\dm_\infty^2\dm\n \overline w^n(t)\dm_2^2\leq c\dm \mathbb P_2\Delta \overline w^n(t)\dm_2\dm \n\overline w^n(t)\dm_2^3\,,\VS \dm\mathbb P_2[\overline w^n\cdot\n v](t)\dm_2^2&\hskip-0.2cm\leq\dm \overline w^n(t)\dm_6^2\dm\n v(t)\dm_3^2\leq c_1\dm\n v(t)\dm_3^2\dm\n\overline w^n(t)\dm_2^2 \leq \frac{c}{t}\dm\n\overline w^n(t)\dm_2^2\,,\VS \dm \mathbb P_2[v\cdot\n \overline w^n](t)\dm_2^2&\hskip-0.2cm\leq \dm v(t)\dm_\infty^2\dm \n \overline w^n(t)\dm_2^2\leq \frac{c}{t}\dm\n\overline w^n(t)\dm_2^2\,.
\ea
\]
Therefore, we arrive at
\[
\frac{d}{dt}\dm \n\overline w^n(t)\dm_2^2\leq c\bigg(1+\frac{1}{s}\bigg)\dm\n\overline w^n(t)\dm_2^2+c\dm \n\overline w^n(t)\dm_2^6\,.
\]
This last inequality can be integrated over a maximal interval of time $(s,s+T(s))$, leading to 
\[
\overline w^n\in C(s,s+T(s);J^{1,2}(\OO))\cap L^2(s,s+T(s);W^{2,2}(\OO)),\,\overline w^n_t,\n\pi_{\overline w^n}\in L^2(s,s+T(s);L^2(\OO)),
\]
with a uniform bound in $n\in\N$. \par Following the procedure employed in the proof of Lemma\,\ref{le: Ex-pert-start}, we deduce the existence of a weak solution $\overline w$ to problem \eqref{pert-start} with initial datum $w(s)$. Moreover, employing a uniqueness argument for the weak limit, we also have
\be\label{reg-LS}
\overline w\in C(s,s+T(s);J^{1,2}(\OO))\cap L^2(s,s+T(s);W^{2,2}(\OO)),\,\overline w_t\in L^2(s,s+T(s);L^2(\OO)).
\ee
Therefore, by a straightforward application of the Prodi-Serrin uniqueness criterion, it follows that $w=\overline w$ on $(s,s+T(s))$. The lemma is completely proved.
\ep
We are now in a position to conclude the proof of Theorem\,\ref{thm: L2-WPpert}, hence we prove a structure theorem for the weak solution $w$, established in Lemma\,\ref{le: Ex-pert-start}.
\begin{lemma}\label{le: Struct-start}
    {\sl Let $w$ be the solution established in Lemma\,\ref{le: Ex-pert-start}. Then, for all $\eta>0$ there exists $\theta\in[T_0,T_0+\eta^{-2}(1-A)^{-1}\dm \mathbb P_2w_0\dm_2^4]$, and a sequence of open time intervals $\{(\theta_\ell,\Theta_\ell)\}_{\ell\in\N_0}$, with $(\theta_0,\Theta_0)=(\theta,\infty)$, such that:
    \[
    \ba{c}
    (0,\infty)\setminus \underset{\ell\in\N_0}\cup (\theta_\ell,\Theta_\ell)\,\,\text{has zero Lebesgue measure, }\VS w\,\text{is a regular solution on }(\theta_\ell,\Theta_\ell)\times\OO\,\text{for all }\ell\in\N_0\,.
    \ea
    \] }
\end{lemma}
        \bp
          We infer that for all $\eta>0$ there exists $\theta\in [T_0, T_0+\eta^{-2}(1-A)^{-1}\dm\mathbb P_2w_0\dm_2^4]\equiv[T_0,\overline T]$ such that
        \[
        \dm \mathbb P_2w_0\dm_2\dm \n w(\theta)\dm_2\leq \eta\,.
        \]
        Assume that the statement is false. Then we have, for all $t\in [T_0, \overline T]$, $\dm \mathbb P_2w_0\dm_2\dm \n w(t)\dm_2> \eta$. Hence, recalling the energy inequality, we deduce
        \[
        \ba{l}\displ
        \dm\mathbb P_2w_0\dm_2^2=(1-A)(1-A)^{-1}\eta^{-2}\eta^2\dm\mathbb P_2w_0\dm_2^4\dm\mathbb P_2w_0\dm_2^{-2}=(1-A)\eta^2\dm\mathbb P_2w_0\dm_2^{-2}\int_{T_0}^{T_0+\eta^{-2}(1-A)^{-1}\dm\mathbb P_2w_0\dm_2^4}dt\VS\hskip5cm <2(1-A)\int_{T_0}^{T_0+\eta^{-2}(1-A)^{-1}\dm\mathbb P_2w_0\dm_2^4}\dm\n w(t)\dm_2^2\,dt\leq \dm\mathbb P_2w_0\dm_2^2\,,
        \ea
        \]
        that is a contradiction. \par We take $w(\theta)$ as initial datum for problem \eqref{pert-start}. By virtue of Lemma\,\ref{le: Uni-LS}, $w$ is a regular solution on some maximal interval of time $(\theta,\theta+T(\theta))$. We show that $T(\theta)=\infty$. For all $t\in (\theta,\theta+T(\theta))$, it holds
        \[
        \frac{d}{dt}\dm \n w(t)\dm_2^2+\frac12\dm \mathbb P_2\Delta w(t)\dm_2^2+\dm  w_t(t)\dm_2^2\leq c\dm \n w(t)\dm_2^6+c(\dm\n v(t)\dm_3^2+\dm v(t)\dm_\infty^2+|u_\infty|^2)\dm \n w(t)\dm_2^2\,.
        \]
        Setting $y(t):=\dm \n w(t)\dm_2^2$ and $a(t):=c(\dm v(t)\dm_\infty^2+|u_\infty|^2)$, we have
        \[
        y'(t)\leq a(t)y(t)+cy^3(t)\,.
        \]
        Since $y(\theta)\leq \frac{\eta^2}{\dm\mathbb P_2 w_0\dm_2^2}$, with $\eta$ arbitrarily small, and, by virtue of Theorem\,\ref{thm: GWP-aux}, the function $a(t)$ is bounded for $t\geq \theta$, employing Lemma\,\ref{le: Mar-AB} we conclude that $T(\theta)=\infty$. \par We now discuss the regularity of the solution on $(0,\theta)$. Let $t_\alpha\in (0,\theta)\cap \widetilde{\mathcal T}_\theta$. By virtue of Lemma\,\ref{le: Uni-LS}, we deduce that $w$ is regular on a maximal interval of time $I_\alpha:=[t_\alpha,t_\alpha+T(t_\alpha))$. Let $t_\beta\in \big((0,\theta)\cap\widetilde{\mathcal T}_\theta\big)\setminus I_\alpha $. Lemma\,\ref{le: Uni-LS} ensures that $w$ is regular on a maximal interval of time $I_\beta:=[t_\beta,t_\beta+T(t_\beta))$. We have two possibilities
        \begin{description}
            \item[i.] $I_\alpha\cap  I_\beta=\emptyset$;
            \item[ii.] $I_\alpha\cap I_\beta\ne \emptyset$, then $I_\alpha\subset I_\beta$.
        \end{description}
        We justify the second item. Since $t_\beta\notin I_\alpha$, then $t_\beta<t_\alpha$. Therefore, we necessarily have $t_\alpha\in I_\beta$ and, in view of the maximality of $I_\beta$, $I_\alpha\subset I_\beta$. \par If item $\bf i.$ holds, we take $t_\gamma\in \big((0,\theta)\cap\widetilde{\mathcal T}_\theta\big)\setminus (I_\alpha\cup I_\beta)$, otherwise  we take $t_\gamma\in \big((0,\theta)\cap\widetilde{\mathcal T}_\theta\big)\setminus I_\beta$. Iterating the procedure, we construct a family of intervals $\{I_\alpha\}_{\alpha\in\mathbb A}$ such that
        \[
        (0,\theta)\cap\widetilde{\mathcal T}_\theta\subseteq \underset{\alpha \in \mathbb A}\cup I_\alpha\,.
        \]
        Since the interior sets of the quoted intervals  are disjoint, then $\mathbb A$ is at most countable. The lemma is completely proved.
\ep
We are now in a position to prove the main result of this work, Theorem\,\ref{thm: MR}.
\bp[Proof of Theorem\,\ref{thm: MR}]
Let $u_0\in L^p(\OO)$, $p\in (2,3)$. Recalling that $\mathbb P_3:L^3(\OO)\to J^3(\OO)$ is a bounded linear operator, we have
\[
\dm \mathbb P_3 g\dm_3\leq c(3)g\,,\,\,\text{for all }g\in L^3(\OO)\,.
\]
In view of Lemma\,\ref{Dec-Lp}, for all $\rho>0$ there exist $v_0\in L^3(\OO)$ and $w_0\in L^2(\OO)$ such that
\[
u_0=v_0+w_0\,;\,\,\dm v_0\dm_3\leq \rho^{1-\frac{p}{3}}\dm u_0\dm_p^{\frac{p}{3}}\,,\,\,\dm w_0\dm_2\leq \rho^{1-\frac{p}{2}}\dm u_0\dm_p^{\frac{p}{2}}\,.
\]
We consider $v_0$ and $w_0$ as initial data for problem \eqref{aux-start} and problem \eqref{pert-start} respectively.\par
Let $\delta$ be the positive constant stated in Theorem\,\ref{thm: GWP-aux} and choose
\be \label{rho}
\rho\leq \bigg(\frac{\delta-|u_\infty|}{c(3)\dm u_0\dm_p^{\frac{p}{3}}}\bigg)^{\frac{3}{3-p}}\,.
\ee
Then, we obtain a weak solution $u$ to problem \eqref{NS-Ose} as the sum $u=v+w$, where, by virtue of Theorem\,\ref{thm: GWP-aux}, $v$ is a regular solution to \eqref{aux-start} and, by virtue of Lemma\,\ref{le: Ex-pert-start}, $w$ is a weak solution to \eqref{pert-start}. It is easily seen that the regularity properties deduced for $v$ and $w$ ensure that $u$ fulfills all the items given in Definition\,\ref{def: WS-Ose}. \par Theorem\,\ref{thm: GWP-aux} and Lemma\,\ref{le: Struct-start} allow us to say that, for all $T>\theta$,
\begin{equation}\label{cond-theta}
    u-U\in C([\theta,T);J^3(\OO))\,.
\end{equation}
Let us consider the following initial boundary value problem:
\begin{equation} \label{pert-staz-PMR}
    \begin{cases}
       \overline u_t-\Delta \overline u+\overline u\cdot\n\overline u+u_\infty\cdot\n\overline u+U\cdot\n\overline u+\overline u\cdot\n U+\n\pi_{\overline u}=0\,,\,\,\text{in }(\theta,T)\times\OO\,,\\ \hskip3.99cm\n\cdot\overline u=0\,,\,\,\text{in }(\theta,T)\times\OO\,,\\ \hskip3.75cm\overline u(t,x)=0\,,\,\,\text{on }(\theta,T)\times\partial\OO\,, \\ \hskip2.8cm\displ\lim_{|x|\to\infty }\overline u(t,x)=0\,,\\ \hskip2.8cm\overline u(\theta)=u(\theta)-U\,,\,\,\,\text{in }\OO\,.
    \end{cases}
\end{equation}
We recall that, in view of Lemma\,\ref{le: WP-staz}, Theorem\,\ref{thm: GWP-aux} and Lemma\,\ref{le: Struct-start}, $\dm u(\theta)-U\dm_3$ satisfies a suitable smallness condition. Hence, by virtue of Lemma\,\ref{le: WP-Shibata}, problem \eqref{pert-staz-PMR} has a unique mild solution $\overline u$ in the class $C([\theta,\infty);J^3(\OO))$ satisfying 
\[
\dm \overline u(t)\dm_q\leq ct^{-\frac32(\frac13-\frac1q)}\,,\,\,\text{for all }3\leq q\leq \infty\,\,\text{and }t>\theta\,.
\]
Since $u-U$ solves \eqref{pert-staz-PMR} and belongs to $C([\theta,\infty);J^3(\OO))$, we easily deduce that $\overline u=u-U$ and we conclude the proof. 
\ep
\begin{rem}
    {\rm We point out that our techniques can be adapted to the case of an initial perturbation $u_0\in J^3(\OO)$. As a matter of fact, given $u_0\in J^3(\OO)$, for all $\vep>0$ there exists $w_0\in \mathscr C_0(\OO)$ such that
    \[
    \dm u_0-w_0\dm_3<\vep\,.
    \]
    Setting $v_0:=u_0-w_0$, we obtain a weak solution to \eqref{NS-Ose} as the sum $u=v+w$, where $v$ and $w$ are constructed as in subsections\,\ref{sec: auxprob}-\ref{sec: pertprob}. The quoted solution has the asymptotic properties given in Theorem\,\ref{thm: MR}. \par Therefore, by our techniques we are able to obtain an analogous result to that found by Hishida and Maremonti in \cite{His-Mar}.}
\end{rem}

{\bf Acknowledgements. }{\small I would like to thank Prof. P. Maremonti for his careful supervision of this work and for suggesting me to look for a cut-off function as that employed in the proof of Lemma\,\ref{le: Conv-extball}.}

\end{document}